\documentclass[10pt,twoside]{article}

\usepackage{amsfonts} 
\usepackage{units, amsmath, color}
\usepackage{pstricks,pst-node,pst-coil,pst-plot}  

\usepackage{graphicx}

\usepackage{amssymb}
\usepackage{latexsym}
\usepackage{hyperref}

\makeatletter

\@addtoreset{equation}{section} \makeatother
\newtheorem{theorem}{Theorem}[section]
\newtheorem{remark}[theorem]{Remark}
\newtheorem{lemma}[theorem]{Lemma}
\newtheorem{proposition}[theorem]{Proposition}
\newtheorem{corollary}[theorem]{Corollary}
\newtheorem{definition}[theorem]{Definition}
\newtheorem{example}[theorem]{Example}

\def\1{\mathbf{1}}
\def\:{\lrcorner}
\def\#{\sharp}

\def\a{\alpha}

\def\d{\delta}

\def\o{\circ}
\def\s{\sigma}

\def\qed{\ensuremath{\quad\Box\quad}}

\def\inv#1{\raise.1em\hbox to 0pt{$^{-1}$\hss}_{#1}\;}

\def\v{\noindent}

\newcommand{\bean}{\begin{eqnarray*}}
\newcommand{\eean}{\end{eqnarray*}}
\newcommand{\benu}{\begin{enumerate}}
\newcommand{\eenu}{\end{enumerate}}
\newcommand{\eea}{\end{eqnarray}}
\newcommand{\bea}{\begin{eqnarray}}

\def \beit{\begin{itemize}}
\def \eeit{\end{itemize}}
\def \bey{\begin{eqnarray*}}
\def \eey{\end{eqnarray*}}

\def \bui#1#2{\mathrel{\mathop{\kern 0pt#1}\limits^{#2}}}
\def \buil#1#2{\mathrel{\mathop{\kern 0pt#1}\limits_{#2}}}

\newtheorem{Theorem}{Theorem}

\newtheorem{Definition}{Definition}

\newtheorem{Corollary}{Corollary}

\newcommand{\be}{\begin{equation}}
\newcommand{\ee}{\end{equation}}

\newcommand{\N}{{\mathbb N}}

\newcommand{\R}{{\mathbb R}}

\newcommand{\SSS}{{\mathbb S}}

\newcommand{\ben}{\begin{enumerate}}
\newcommand{\een}{\end{enumerate}}
\newcommand{\bit}{\begin{itemize}}
\newcommand{\eit}{\end{itemize}}
\newcommand{\edoc}{\end{document}}

\newcommand{\bdefi}{\begin{definition}}
\newcommand{\btheo}{\begin{theorem}}
\newcommand{\bprop}{\begin{proposition}}
\newcommand{\brema}{\begin{remark}}
\newcommand{\bcoro}{\begin{corollary}}
\newcommand{\blemm}{\begin{lemma}}
\newcommand{\bexam}{\begin{example}}

\newcommand{\edefi}{\end{definition}}
\newcommand{\etheo}{\end{theorem}}
\newcommand{\eprop}{\end{proposition}}
\newcommand{\erema}{\end{remark}}
\newcommand{\ecoro}{\end{corollary}}
\newcommand{\elemm}{\end{lemma}}
\newcommand{\eexam}{\end{example}}

\title{Which spacetimes admit conformal compactifications?}

\begin{document}

\author{Olaf M\"uller\footnote{Fakult\"at f\"ur Mathematik, Universit\"at Regensburg, D-93040 Regensburg, \texttt{Email: olaf.mueller@ur.de}}}

\date{\today}
\maketitle

\begin{abstract}
\v We consider the future causal boundary as a tool to find obstructions to conformal extensions, the latter being a slight generalization to conformal compactifications. 
\end{abstract}

\begin{small}{\it Mathematics Subject Classification} (2010):  53A30, 53C50, 53C80\\
\end{small}

\section{Introduction}

\v In mathematical relativity, a question often considered is the existence of a {\em $C^k$-conformal compactification} of a given spacetime $(M,g)$, that is, of an open conformal embedding of the spacetime into a bigger one with a Lorentzian metric of class $C^k$ whose image is relatively compact and satisfies varying additional causal assumptions. One important purpose of conformal extensions is the definition of future null infinity. Maybe the most famous example of a smooth conformal compactification is the {\em Penrose embedding} of Minkowski spacetime into the Einstein cylinder arising naturally from the stereographic projection (see \cite{BEE}, e.g.). This compactification has been used several times in the analysis of certain conformally covariant PDEs. Choquet-Bruhat and Christodoulou \cite{CC}, e.g., used the Penrose embedding to prove existence of global solutions for Yang-Mills-Higgs-Dirac systems in Minkowski space for initial values small in a weighted Sobolev norm. If one tries to generalize this technique to other relevant spacetimes one encounters the difficulty that asymptotically Schwarzschildian relativistic initial values with nonvanishing mass obstruct the existence of a conformal compactification due to the existence of a conformal singularity at spatial infinity (where the Weyl curvature diverges). Therefore, in a recent article of the author with Nicolas Ginoux \cite{GM}, the notion of '$C^k$-conformal compactification' has been slightly generalized to the notion of '$C^k$-conformal extension' in order to accommodate nonzero mass as well, without ever assuming some kind of asymptotical flatness in the sense of curvature tending to zero. Doing so, it was possible to obtain statements similar to those of Choquet-Bruhat and Christodoulou for maximal vacuum Cauchy developments of asymptotically Schwarzschildian data instead of Minkowski spacetime. Thus let us review the definition following \cite{GM}: A subset $A$ of a spacetime $M$ is called {\bf future compact} iff it is contained in the past of a compact subset of $M$, and it is called {\bf causally convex} iff there is no causal curve in $M$ leaving and re-entering $A$. A {\bf conformal extension} $E$ of a globally hyperbolic spacetime $(M,g)$ is an open conformal embedding of $(I^+(S_0),g)$ (where $S_0$ is a Cauchy surface of $(M,g)$) into another globally hyperbolic spacetime $(N,h)$ with a Lorentzian metric of class $C^k$ such that the closure of the image is future compact and causally convex. This generalizes the usual notion of 'conformal compactification' inasfar as it solely requires {\em future} compactness of the closure of the image instead of compactness. 

\bigskip

\v From the work of Friedrich, Anderson-Chrusciel, Lindblad-Rodnianski and others (\cite{hF}, \cite{AC}, \cite{D}, \cite{DS}, \cite{LR}, \cite{SB}) it follows that there is a weighted Sobolev neighborhood $U$ around zero in the set of vacuum Einstein initial values such that, for any $u \in U$, the maximal vacuum Cauchy development of $u$ admits a conformal extension\footnote{In the corresponding theorem in \cite{AC}, the existence of a conformal extension of the spacetime and its time-dual are shown, which in general do not yield a conformal compactification due to the conformal singularity at $i_0$ for nonzero mass because of diverging Weyl curvature.}.

 \bigskip
 
\v Now, given a conformal extension, it is easy to define its future boundary. One can see that this conformal future boundary is homeomorphic to the intrinsic future boundary defined as in \cite{FHS1}, e.g. (and, consequently, its homeomorphy class is independent of the choice of conformal extension). Thus it seems natural to look for obstructions against conformal extensions by considering the intrinsic boundary. As conformal extendibility is weaker than conformal compactifiability, all these obstructions obstruct conformal compactifiability as well.
 
We will obtain two main results about nonexistence of conformal extensions for standard static manifolds $(\R \times S, - dt^2 + k)$. The main tool to obtain the results is a particularly simple metric on the set of indecomposable past sets (IPs) as introduced by Geroch, Kronheimer and Penrose \cite{GKP} that will be defined in the second section. The first of the results, obtained in the third section, requires an additional topological assumption on $S$:
 
 \begin{Theorem}
\label{topological}
Let $M$ be a standard static globally hyperbolic manifold whose Cauchy surface is not homeomorphic to a cone (satisfied e.g. if it is of dimension at least two and has more than one end). Then $M$ has no conformal extension.
\end{Theorem}
 
The second main result, obtained in the fourth section, requires an additional geometric assumption on the metric $k$ of the standard static slice:
 
 \begin{Theorem}
 \label{static}
A standard static spacetime over a complete Riemannian manifold $(S,k)$ can only have a conformal extension if the Busemann boundary of $(S,k) $ is the Gromov boundary of $(S,k)$.
\end{Theorem}

 \bigskip
 
\v It is interesting to compare $C^k$-conformal extendibility with Penrose's definition of $k$-asymptotic simplicity (Def. 9.6.11 in \cite{PR}). First, $C^k$-conformal extendibility 'almost' implies $k$-asymptotic simplicity: Adding to $M$ the future boundary $J^+$ of a conformal extension and its time-dual $J^-$ yields a (Lipschitz) manifold-with-boundary $J^+ \cup  J^-$ satisfying the assumptions of Penrose's definition except for the condition that the gradient of the conformal factor vanishes at $J^{\pm}$. On the other hand, asymptotic simplicity is to weak for the analytic applications in \cite{GM}. However, it is interesting that for $J^{\pm}$ being null, one obtains a cone structure on $J^{\pm}$ as well (see Theorem 9.6.19 in \cite{PR}), just as in the case of a conformal extension of a standard static globally hyperbolic manifold (cf. the proof of Theorem \ref{topological}).

\bigskip

The results on conformal extensions of synoptic globally hyperbolic manifolds obtained in the second section have interesting implications for the maximal vacuum Cauchy developments of initial values close to the trivial ones. For the construction of Hadamard states in the context of CCR quantization, it is very useful that the extension be {\bf essentially null} in the sense that the boundary of the image is a smooth null hypersurface except for a set of measure zero. In that case, the usual symplectic form $\omega$ can, for linear theories, be expressed by Stoke's Theorem as an integral along the future boundary $J^+$ as follows: $\omega (f_1,f_2) = \int_{J^+} G(f_1) \partial_u G(f_2) - G(f_2) \partial_u G(f_1)  $ (where $u$ is the local null direction of the boundary and $G$ is the Green's function)\cite{vM}[Theorem 4.1]\footnote{Note that the author of this article operates with a stronger definition of 'future boundary' implying in particular smoothness of the boundary. An inspection of the proof, however, shows that it is enough to invoke Stokes' theorem for Lipschitz boundaries \cite{wML}[Theorem 3.34] to make the proof work for the case of an essentially null boundary.}, and this representation is of fundamental importance in several applications in CCR quantization (see, for example, \cite{DHP}). Now, one of the corollaries of Section 2 is the following, givig a partial answer to the question of existence of essentially null conformal extensions:

\begin{Corollary}
\label{Hadamard}
There is a Sobolev-open neighborhood $U$ of the trivial Einstein vacuum initial values such that, for each $u \in U$, the maximal vacuum Cauchy development $M(u)$ has an essentially null conformal extension $f$.
\end{Corollary}

In the last section we draw a simple consequence for FLRW (Friedmann-Lema\^itre-Robertson-Walker) spacetimes (which is interesting to compare to the result in \cite{DHP}):

\begin{Corollary}
\label{cor}  
A globally hyperbolic FLRW spacetime $(M,h):= (\R \times Q, - dt^2 + a^2(t) g)$ has no essentially null conformal compactification if $\int_0^\infty a(t) dt < \infty$ or $\int_{- \infty}^0 a(t) < \infty$. If $\int_0^\infty a(t) dt , \int_{- \infty}^0 a(t) = \infty$, then $(M,h)$ has no conformal extension if one of the following holds:
\begin{enumerate}
\item{$Q$ is not homeomorphic to a cone,}
\item{The Busemann boundary of $(Q,g)$ does not coincide with its Gromov boundary.}
\end{enumerate} 
\end{Corollary}

\v  For an interesting and completely different approach to a conformal extension problem in the context of Cartan geometries, cf. \cite{cF}, where e.g. some quotients of Anti-deSitter space are shown to admit no nonsurjective open conformal embedding at all (Cor. 5.1 of \cite{cF}).

\section{Simple and useful metrics on the set of IPs}

\v First of all we clarify what we consider the future part of the boundary:

\begin{Definition}
The {\bf future boundary} $\partial^+ A $ of a subset $A$ of a spacetime $X$ is defined as the set $\{ x \in \overline{A} \setminus A \vert I^- (x) \cap A \neq \emptyset \} $. 
\end{Definition}

\v Now it is not surprising that if $A$ is open, the causal properties of $A$ and $\partial^+A$ are strongly related:

\begin{Theorem}
If $A$ is open and precompact, then it is causally convex if and only if $\partial^+ A$ is achronal.
\end{Theorem}

\v {\bf Proof.} The implication from the left to the right in the first statement is well-known, for the reverse implication assume that there is a future causal curve $c: [0, T] \rightarrow X$ leaving (at a point $a = c(s) $ with $s>0$) and re-entering $A$. Any maximal $C^0$ future timelike extension $k$ of $c$ will leave $A$ a second time (at $b = k(u)$, say). Now take a convex open neighborhood $V$ of $a$ and an open neighborhood $U$ of $c(0)$ in $A \cap V$ as well as some point $c (t)$ with $t>s$ within $V$. Then there is a {\em timelike} and geodesic future curve $\gamma$ from a point $x$ in $U$ to $c(t)$ intersecting the future boundary at another point $w$ in $V$. Now, the curve $h:= k\vert_{[t, u ]} \o \gamma$ is future causal, somewhere timelike and joins $w$ with $b$. By the push-up lemma \cite{pC} we can find a timelike future curve from $w$ to $b$, contradiction (note that as well $w$ as $b$ are in the {\em future} boundary of $A$).  \hfill \qed

\bigskip

\v In the case of $A$ being the image of a conformal extension, the future boundary is completely prescribed:

\begin{Theorem}
\label{FutureEqualsCauchy}
Let $f: (M,g) \rightarrow (N,h)$ be a conformal extension, then $\partial^+ (f(M))$ is homeomorphic to a Cauchy surface of $(I^+(S_0),g)$. 
\end{Theorem}

\v {\bf Proof.} As the closure of the image of $f$ is causally convex, $\partial^+ (f(M))$ is achronal, and for any Cauchy surface $S$ of $(I^+(S_0), g)$, $f(S)$ is acausal by openness and causal convexity, and as the closure of the image is compact, the closure of $f(S)$ in the ambient space $N$ is compact and again acausal. By \cite{BS3}[Prop. 3.6] there is a topological Cauchy surface $\Sigma$ of $(N,h)$ containing $f(S)$ and a $C^0$ Cauchy time function $t$ with $\Sigma = t^{-1} (\{0\})$. Then there is a continuous function $F \in C^{\infty} (\Sigma, \R^+)$ such that, in the identification $N = \R \times S$ via $t$ we get $\partial^+ (f(M)) = \{ (F(s) , s) \vert s \in f(S) \} $: e.g., we can define $F(s):= {\rm sup} \{  r \in \R \vert (r, f(s) ) \in f(M) \} $. Now obviously, $(F(s) , s)$ is not only in $\partial (f(M))$ but in $\partial^+ (f(M))$. 
\hfill \qed

\bigskip

\v In particular, for any two conformal extensions of a g.h. manifold $(M,g)$, their future boundaries are homeomorphic.

\bigskip

\v Now let us turn to the intrinsic definition of the future boundary. Let $\hat{M}$ be the set of indecomposable past subsets (IPs) of $M$, where a past subset is called {\bf indecomposable} iff it is not the union of two different nonempty past subsets. Geroch, Kronheimer and Penrose, in their seminal article \cite{GKP}[Th. 2.1], show that each IP is necessarily the past of a timelike future curve $c$. If $c$ is $C^0$-extendible in $M$, its past is just the space of its endpoint $p \in M$, which gives an embedding of $M$ into $\hat{M}$, more on that below. If $c$ is $C^0$-inextendible, it is called {\bf terminal} or TIP, for short. Now one can, of course, equip $\hat{M}$ with the causal relation coming from inclusion, and this causal relation obviously extends that of $M$. How to put an appropriate topology on $\hat{M}$ has been a subject to an extensive debate, and finally found a satisfactory but still quite technical answer \cite{FHS1}.

\v Here we want to introduce a particularly easy definition of a metrizable topology on the set of IPs tailored only for the special case of globally hyperbolic manifolds $M$ but at the same time avoiding some technical subtleties of the established definition for the general case as given in \cite{FHS1}. The so obtained `poor man's' boundary $\partial^+ M$ will be homeomorphic to the future boundary of any conformal extension if the latter exists. This is done as follows:

We use Lemma 3.2 of \cite{CGM} stating that, for a globally hyperbolic manifold $(N,h)$ and for any compactly supported $\psi \in C^0 (N, [0, \infty))$ with $\int_M \psi (x) d {\rm vol} (x) = 1$, the function $\tau_{\psi}$ with $\tau_\psi (p) := \int_{I^-(p)} \psi d{\rm vol}_h$ is continuously differentiable with ${\rm grad} \tau_\psi$ timelike or zero everywhere. We will see below that the gradient is timelike on the interior of the support of $\psi$. Therefore, we choose a countable covering of $F:= I^+(S)$ by open precompact sets $U_i$ and define, for $\phi_i \in C\infty (I^+(S))$ with $ \phi_i^{-1} (0) = M \setminus U_i $, 

$$ \phi = \sum_{i \in \N} 2^{-i} (\vert \vert \phi_i \vert \vert_{C^1} + \vert \vert \phi \vert \vert_{L^1})^{-1} \phi_i  $$ 

It is easy to see that $\phi \in C^1 \cap L^1 $ with $\vert \vert \phi \vert \vert_{C^1(F)} , \vert \vert \phi \vert \vert_{L^1(F)} \leq 1$, and that the gradient of $\phi$ is timelike everywhere. We rescale $\phi$ such that $\vert \vert \phi \vert \vert_{L^1(F)} =1 $. 

We can induce a metric on $\hat{M}$ by 

$$\delta(A,B) :=  - {\rm ln}(1- \int_{\Delta(A,B)} \phi dvol ) $$ 

where $\Delta(A,B) := (A \setminus B) \cup (B \setminus A)$ is the symmetric difference of $A$ and $B$. A second, strictly larger, metric $d$ on $\hat{M}$ is defined by 

$$d(A,B) :=  - {\rm ln} (1 - \int_{\Delta(A,B)} \phi dvol ) + {\rm ln} (\int_{\Delta(I_j^+ (A), I_j^+ (B))} \phi d vol)$$

\v where, for a subset $ Q $ of $M$, $I_j^+ (Q) := \cap_{q \in Q} I^+ (q) $ denotes the joint future of $Q$. Note that $ I_j^+ (I^-(p))  = I^+(p)$ for any $p \in M$.

\bigskip

Both metrics are obviously nonnegative and symmetric. The triangle inequality follows immediately from the set-theoretic triangle inequality $\Delta(A,C) \subset \Delta(A,B) \cup \Delta (B,C) $ and convexity and monotonicity of the logarithm. It remains to show that $\d$ does not vanish between different subsets. 

\bigskip

As a corollary of the Geroch-Kronheimer-Penrose structural result, we have that IPs are always open subsets. But in general, for $A$,$B$ two different open past subsets, we have $(A \setminus \overline{B}) \cup (B \setminus \overline{A}) $ is nonempty (and open): Let, w.r.o.g., be $x \in A \setminus B$, then, as $A$ is open,  there is a $y>x$ with $y \in A$, and by the push-up lemma \cite{pC} and the fact that $B$ is past we know that $y \notin \overline{B}$. This shows that the metrics are not pseudometrics but true ones.

\bigskip

A function is called {\bf time function} iff it is monotonically increasing along future timelike curves and it is called {\bf temporal} iff it is $C^1$ and has timelike past gradient. It is called {\bf Cauchy} if its restriction on any $C^0$ inextendible timelike future curve is surjective onto the real numbers. 

One nice feature of $\d$ and $d$ is that we can construct a Cauchy time function which is basically the distance to $i^+$. For $h: (0, \infty) \rightarrow \R $ with ${\rm lim}_{x \rightarrow 0} h(x) = \infty$ and ${\rm lim}_{x \rightarrow \infty} h(x) = - \infty$, we simply define  

$$ t (A):=  h(d(A, M )), \qquad T(A) := h(\d(A,M) ). $$

\begin{Theorem}
$t$ and $T$ are temporal functions. $t$ is always Cauchy, whereas $T$ is never Cauchy.
\end{Theorem}

{\bf Remark:} The relevance of $T$ in spite of not being Cauchy lies in its property that it does not diverge towards large parts of the boundary, therefore it can be taken as some intrinsic replacement for a temporal function of the exterior space along null future infinity. 

\v {\bf Proof.} $t$ is almost the $C^1$ Cauchy temporal function constructed in \cite{CGM}. The only differences are that in \cite{CGM} we have $\vert \vert \phi \vert \vert = \infty$ and $h: \R \rightarrow \R$, $h(x) = -x$, but these differences do not affect the proof, which can be adopted verbatim. The statement on $T$ follows from the proof of the results in \cite{CGM}, in particular Theorem 1.1, Cor. 5.5 and Corollary 5.6. To show that the derivative of $T$ does not vanish basically boils down to the computation of the diamond volume at the tip in terms of $\dot{g}$ and $a$ in a metrical splitting $-a \cdot d \theta^2 + g_\theta$ for an auxiliary smooth Cauchy temporal function $\theta$. \hfill \qed



\bigskip

\v The metrics are, of course, not natural, as the choice of $\phi$ is highly arbitrary. Now, as mentioned above, $I^+(S)$ is embedded in $H$ via a map $s: p \mapsto A:= I^- (c) $ where $c$ is any timelike future curve from $S$ to $p$. The map $s$ maps into the IPs $A$ such that $A \cap S$ is compact, but in general it is not surjective on them: Consider Schwarzschild or Kruskal spacetimes where there are future precompact IPs $A$ that are pasts of inextendible curves and not of points on the spacetime. The map $s$ is injective because $(M,g)$ is distinguishing.

 Now we define $ \partial^+ M := \hat{M} \setminus I(M) $, equipped with the subspace topology.  

We also define the map $E$ associating to a shadow $I^-(c)$ the endpoint of $f \o c$ in $N$. We denote the curve containing the endpoint as $\overline{c}$. The choice of $c$ is of course not canonical, but for two curves $c,k$ with $ I^-(c) = I^-(k) $ we have 

\bean
I^-(\overline{k}) = I^- (k ) = I^- (c) = I^- (\overline{c})
\eean 

and thus $E(c) = E(k)$ as $(N,h)$ is distinguishing. Thus $E$ is well-defined. The next theorem shows that the intrinsic future boundary defined above is homeomorphic to the future boundary of a conformal extension. For the sake of its proof, we recall the notion of inner and outer continuity of set-valued maps: a set-valued map $F$ is called {\bf inner (resp., outer) continuous at a point $p$} iff for all compact sets $C \subset {\rm int}(F(p)) $ (resp., for all compact sets $ K \subset M \setminus \overline{F(p)}$), there is an open set $W$ containing $p$ such that for all $q \in W$, we have $C \subset {\rm int}(F(q))  $ (resp., $K \subset M \setminus \overline{F(q)}$).

\begin{Theorem}[compare with Th. 4.16 from \cite{FHS1}]
\label{hadescauchy}
Let $f:(M,g) \rightarrow (N,h)$ be a conformal extension. Then $E$ is a homeomorphism between $\partial^+ M$ and $\partial^+ (f( M)) $. Its inverse is the map $\overline{s}: p \mapsto f^{-1} (I_N^-(p))$. 
\end{Theorem}

\v {\bf Proof.} We want to show that $\overline{s}$ is a right and left inverse of $E$. First we have to show that $\overline{s}$ takes values in the IPs, but this is just a consequence of omitting the final point and using the openness of timelike future cones in the ambient manifold, thus $f^{-1} (I_N^- (\overline{f \o c})) = f^{-1} (I_N^- (f \o c)) = I^-_M (c)$.  Therefore, indeed, for $q \in \partial^+ f(M)$, the set $f^{-1} (I^-(q) \cap f(M) ) $ is an IP in $M$. And $\overline{s}$ is a right inverse of $E$ as $N$ is distinguishing, it is a left inverse of $E$ as $c$ generating the IP is a causal curve.

The map $\overline{s}$ is continuous: First, the assignment $ p \mapsto I^- (p) $ is inner continuous in any spacetime and is outer continuous in causally continuous spacetimes and globally hyperbolic manifolds are causally continuous (see \cite{BEE}, e.g.).  Then, the inner and outer continuity and the fact that the image consists of precompact sets imply continuity in $d$, as $\overline{B}_r ( C)$ and $M \setminus B_r (M \setminus C)$ are compact sets for a precompact set $C$ and for any $r >0$, and as $\int_M \phi = {\rm lim}_{n \rightarrow \infty} \int_{K_n} \phi $.  \hfill \qed

\bigskip

Now we want to examine the structure of the boundary a bit closer in the case of {\em synoptic} manifolds. Here, in accordance with \cite{oM3}, a spacetime $(M,g)$ is called {\bf synoptic} iff the timelike futures of any two points intersect each other. This property is equivalent to the existence of a timelike future curve $c$ with $I^-(c) = M$. In the usual terminology of the causal boundary this property is often also called {\em indecomposable} ($M$ considered as a past subset of itself). 

\begin{Theorem}
\label{timelike-infinity}
Let $f$ be the conformal extension of a synoptic globally hyperbolic manifold $(M,g)$ in $(N,h)$. Then there is a point $i_+ \in N$ such that for all Cauchy surfaces $S$ of $I^+(S_0)$ we have $f(I^+_M(S)) = I^- (i_+) \cap I_N^+ (f(S)) $. 
\end{Theorem}

\v {\bf Proof.} Define $i^+ := E(f(c))$ for any maximal curve $c$. Then one inclusion is obvious, the other is a direct consequence of $f(M) $ being causally convex. \hfill \qed  

\bigskip

\v Obviously, $ T$ as defined above does only diverge towards $i^+$ and does not diverge towards the rest of the boundary. Thus it can serve as an analogon of a Cauchy temporal function of the surrounding nonphysical space. The Cauchy temporal function $t$ of course diverges towards any point of the boundary.

 
 \bigskip

\v Now we can give the announced proof of the first corollary:

\bigskip

{\bf Proof of Corollary \ref{Hadamard}.}
 The stability result of Lindblad-Rodnianski \cite{LR} states that there is a neighborhood $U$ of the trivial initial values open in some Sobolev norm such that, for all $u \in U$, the maximal vacuum Cauchy development of $u$ is timelike geodesically complete, and gives concrete estimates for them that imply good convergence to the Minkowski metric. Inspection of their formulas immediately gives that the past of any $x_0$-coordinate line is all of $M$, thus the developments are synoptic. It is known that, restricting further to a smaller open subset if necessary, all of them have conformal extensions, as mentioned in the introduction. Theorem \ref{timelike-infinity} implies that the boundary of $I^-(i^+) $ in $I^+ (S) \subset N $. But then the boundary of $I^- (i^+)$ is essentially null: Except in the zero measure set given by the cut locus of $i^+$, the boundary at some point $p = exp(n)$ is smooth and contains the null vector $d_n exp \cdot n$.  \hfill \qed 



\bigskip


\section{Topological obstructions against conformal extendibility}



\bigskip

Now, considering Corollary \ref{Hadamard}, one could try to prove that the future boundary of any conformal extension of a synoptic globally hyperbolic manifold is homeomorphic to a cone, restricting thereby the possible topologies of the Cauchy surfaces. This is, however, wrong, as can be seen in the following counterexample: Consider the Lorentzian direct product $M := \R \times \SSS^1 \times \R$ and choose $p \in M$ arbitrarily. Then the past cone $I^-(p)$ is easily seen to be synoptic, conformally extendable by the inclusion and to have a Cauchy surface diffeomorphic to $\SSS^1 \times \R$ (thus not homeomorphic to a cone). We can, however, prove the presence of a cone structure in the case of $(M,g)$ being standard static. In that case, the boundary can be described using the Gromov compactification which we now define.

\v The {\bf Gromov compactification} of a complete metric space $Z$ is the closure of the distance functions plus constants in the function space topology (pointwise convergence or uniform convergence or Lipschitz convergence, all of those are equivalent on the space $L_1 (Z , \R)$ of Lipschitz functions of Lipschitz constant $1$). Explicitly: We have a map $j \rightarrow L_1 (Z, \R)/ \R $ given by $j(x) := [ d^g( x, \cdot) ]$, and we define the Gromov compactification $Z_G$ by $Z_G := \overline{j(Z)}$. It is easy to see that $ L_1 (Z, \R)/ \R $ is a compact metrizable space (see \cite{FHS1}, e.g.). Consequently, $Z_G$ is also a compact metrizable space. Now, assuming that $Z$ has a differentiable structure, for any unit-speed-or-lower piecewise $C^1$ curve $c:[0,D) \rightarrow Z$ (with $ D \in \R \cup \{ \infty \} $), we define the Busemann function $b_c$ of $c$ by $b_c (x) := \lim_{t \rightarrow D} (t - d(x, c(t)))$. It is an easy exercise to show that if this limit is finite for some $x$, it is finite everywhere and defines a Lipschitz function of Lipschitz constant $1$. 

\bigskip

\v In the case of $(M,g)$ being standard static, one can now give another characterization of $\hat{M}$. In particular, in this case the shadows corresponding to true boundary points are precisely the ones whose closure in $S$ is not compact:

\begin{Theorem}
\label{cone}
Let $(M,g) = (\R \times S, -dt^2 + g_0)$ be standard static, then $(\hat{M}, d) \setminus \{ M \}$ is homeomorphic to the subset of equivalence classes of distance functions and Busemann functions on $S$ in the Gromov compactification with the function space topology. \end{Theorem}

\v {\bf Proof.} Each point $(t,q)$ defines a distance function $d^g (q) +t $, and more generally, each future timelike curve $c(s) = (s, q(s))$ associated to an IP defines a Busemann function $b_q$. Conversely, to each Busemann function $f$, we can associate the set $A_f := \{ (s,x) \vert s < f(x) \}$. It is obvious that $I^-(c) = I^-(k) \Leftrightarrow b_c^{-1} (0) = b_k^{-1} (0)$, and each Busemann function $b$ is uniquely fixed by its zero level set $K$, as $b(p) = d(p,K)$. The assignment is continuous for the metric $d$ on $\hat{M}$ as the assignment of sublevel sets to Lipschitz functions is inner and outer continuous. \hfill \qed


\bigskip

\v As an additional structure, there is an $\R$-action $\cdot_t$ on $M$ related to the standard static time function $t$, given as $r \cdot_t (s, x) = (s+r, x)$ in the splitting of $M$ given by $t$. The action can be extended to $\hat{M}$ in three equivalent ways:

\begin{itemize}
\item On the level of IPs as $(r, A) \mapsto r \cdot_t A $,
\item On the level of associated curves as $(r,c) \mapsto r \cdot_t c$,
\item On the level of Busemann functions as the action by adding constants: $(r,b) \mapsto r+b $. 
\end{itemize}

\v The correspondence is easy to see - in the last case it is just given by taking zero-sublevel sets of the Busemann functions. The action is proper and free on the complement of $i^+ = M$ and has a section $\s : B(S) \rightarrow \{  b \in B(S) \vert b(q) = 0 \} $ given by subtracting the respective values at a fixed point $q$, i.e. $\s (b) = b - b(q) $ (see e.g. \cite{FH}). Furthermore, all orbits of the $s$-action approach $i^+$.

\v Consequently, the set of TIPs is homeomorphic to a cone over $ B(S) / \R $, where the $\R$-action is given by addition of constants.

\bigskip

\v {\bf Proof of Theorem \ref{topological}.} First we show the statement in the brackets: If a Cauchy surface $S$ is of dimension at least two and has more than one end, then there is a compact subset $K \subset S$ such that $S \setminus K $ is disconnected. In contrast, in a cone over a {\em connected} topological space there is no such compact subset, whereas the cone over a {\em disconnected} topological space is never a topological manifold. 

\v Now, assume that there is a conformal extension. Then Theorems \ref{FutureEqualsCauchy} and \ref{hadescauchy} imply that $\partial^+M$ is homeomorphic to $S$ and thus not homeomorphic to a cone. On the other hand, Theorem \ref{cone} and the remarks following it imply that $\partial^+M$ is homeomorphic to a cone, contradiction.  \hfill \qed  


\section{A geometric obstruction}

 Now, in many cases, not every point of the Gromov completion of a complete Riemannian manifold is a Busemann function. A famous example discovered by Steven Harris is the unwrapped-grapefruit-on-a-stick \cite{H}. A grapefruit-on-a-stick is just a rotational surface with profile curve equal to $1/10$ outside $[-1,1]$, equal to a semicircle in $[-1 + a, 1 - a]$ for $ a < 1/10$ and interpolating smoothly in the remaining two intervals. Let $M$ be the Riemannian universal cover of this manifold and call it the unwrapped-grapefruit-on-a-stick. The next theorem tells us that such a behavior obstructs the existence of a conformal extension. This is because if a globally hyperbolic manifold has a conformal extension, the future boundary of the image is also future-compact, thus for any compact set $C$ in $S$, $J^+ (C) \cap \partial^+M  $ is compact. Translated into the terminology of the IPs that means that for any point $p \in S$, the subset $H_p$ of all IPs containing $p$ must be compact. We want to exploit this fact in the case of standard static spacetimes using Busemann functions:

\bigskip

\v {\bf Proof of Theorem \ref{static}.}  As the Gromov compactification is the closure of $j(M)$  in the function space topology and as $j(M)$, in turn, is a subset of the Busemann functions, future-compactness of the future boundary (a necessary condition as explained above) is only possible in the case that all points in the Gromov compactification are Busemann functions: Assume a sequence of Busemann functions $b_{c_n}$ (of finite or infinite lengths) converges in the function space topology to a point $F$ in the Gromov compactification, and pick any point $p \in G$, then all $ c_n $ have to take their image eventually in $I^+ (\overline{B_{F(p)+1} (p)}) $, i.e., in the future of a compact set.  \hfill \qed

\bigskip

\v Consequently, we have

\begin{Corollary}
The standard static spacetime over the unwrapped-grapefruit-on-a-stick does not have a conformal extension. 
\end{Corollary}

\v {\bf Proof.} The sequence of functions $n - d((0, 2 \pi n)) $ does not converge to a Busemann function as shown in \cite{H}. We repeat the proof here in more detail: We coordinatize the manifold by the lift of cylindrical coordinates on the grapefruit-on-a-stick. Then consider $ A:= (0, -1), B:= (0,0), C:= (0,1)$. And obviously $ \vert d(A, y) - d(C, y) \vert \geq 2$ for all $y \in M \setminus \R \times (-1/3, 1/3)$. Therefore if the sequence $n - d((2 \pi n,0), \cdot)$ converges to the Busemann function of a curve $k$, then $k$ has to be eventually in $\R \times (-1/3, 1/3)$. But then, curve length increases as the first coordinate but distance to $(0,0)$ only as $1/10$ thereof, thus the limit in the definition of the Busemann function cannot exist. \hfill \qed

\bigskip

\v {\bf Remark:} Note that the unwrapped-grapefruit-on-a-stick has one end and is of bounded geometry, being a Riemannian cover of a compact Riemannian manifold, and one can even modify the construction to get an asymptotically Euclidean example displaying the same behaviour.
 
\bigskip

\v Let us consider an example given in \cite{FH}: Let $M$ be the warped product with a compact manifold: $M = (\a, \omega) \times_a K$, where $\alpha, \omega \in \{ - \infty \} \cup \R \cup \{ \infty \}  $, $a : (\alpha,\omega) \rightarrow \R$ is a positive function, $K$ is compact, and the metric on $M$ is $h = dr^2 + a(r)^2 h_K$ for a Riemannian metric $h_K$ on $K$.
If there is an $s \in (\a, \omega)$ such that $a$ is monotonically decreasing when restricted to $(\a, s)$ and 
monotonically increasing when restricted to $(s, \omega)$. Then Theorem 6.2 from \cite{FH} assures that 
$\partial B(M)$ consists of two spaces $B_{\alpha} $ (attached at $\{ \a \} \times K $) and $ B_{\omega}$ (attached at $\{ \omega \} \times K$), and for $\iota \in \{ \a, \omega \} $ we have: $B_{\iota}$ is homeomorphic to $K$ if and only if $ \vert \int_s^{\iota} 1/a(s)^2 dr \vert < \infty$, otherwise, $B_{\iota}$ is a single point.

\v Here, the Busemann boundary is obviously compact, and thus coincides with the Gromov boundary. If we choose $s= \a$ and assume that $M$ can be extended to a Riemannian manifold $N$ by adding one point at the end $\a$, then $N$ has only one end but still the corresponding standard static manifold does not have a conformal extension as $M$ is not homeomorphic to a cone over a point (i.e., to the real line). Thus, the linear increase of the warping factor in the Euclidean space is not quite a kind of limiting case, but for every exponent greater than $1/2$, the criterion above would still allow for a conformal extension.  

\bigskip

\v Note that the previous results are not of Riemannian nature but genuinely Lorentzian as one can choose metrics on the Cauchy surfaces that have conformal extensions following the results of Marc Herzlich \cite{mH}.

\section{Consequences for FLRW spacetimes and comparison to the usual topology on the set of IPs}

\v Friedmann-Lemaitre-Robertson-Walker spacetimes $(\R \times N, - dt^2 + a(t) g)$ form a class of manifolds often used as cosmological models. 
Now, Corollary \ref{cor} is a simple consequence of the results above for those manifolds:
\bigskip

{\bf Proof of Corollary \ref{cor}}: First of all, $g$ is complete, as the spacetime is globally hyperbolic. Now assume that $\int_0^{\infty} a(t) dt < \infty$. In that case, for all IPs $A$ and every Cauchy surface $S$, the set $A \cap S$ is precompact. On the other hand, given a conformal compactification $i$, choose a Lorentzian normal neighborhood $N$ of $p \in \partial i(S)$. As the compactification is essentially null, there is a point $q = exp_p(n) \in N$ for $n$ null that is contained in $\partial i(M)$, and its past is an IP whose intersection with $S$ is not precompact, contradiction. The other statements are easily proven by performing the obvious conformal change, reparametrizing the $t$-axis by arc length and applying the previous results to the so constructed standard static manifold.    \hfill \qed

\bigskip

\v Finally, it is interesting to compare the future boundary defined as above to the future boundary defined in the framework of \cite{FHS1}, \cite{FHS2} where (for the standard stationary case) the chronological topology is used. In this case, the (forward) Busemann boundary is Hausdorff if and only if it coincides with the Gromov boundary, and in this case one can show (see Theorems 6.10, 6.26 of \cite{FHS2} and Theorem 6 from \cite{H} for the static case) that the chronological topology coincides with the Lipschitz topology, and the future boundary (defined as in the following proof) is homeomorphic to a cone over the Busemann boundary.

\begin{Theorem}
Let $i: M \rightarrow N$ be a conformal extension. Then the future boundary $\partial_i^+ M $ of $i$ is isocausally homeomorphic to the future part $\partial^+ M$ of the causal boundary $\partial M$ of $M$. 
\end{Theorem}

\v {\bf Proof.} As $i(M) $ is precompact, Proposition 4.8 from \cite{FHS1} implies that, in the terminology of that article, $i$ is a chronologically complete envelopment. The isocausal homeomorphism is just the restriction of the map $\pi$ in \cite{FHS1} on the equivalence classes of subsets of the form $(\emptyset, P)$ for a past set $P$. Looking at the corresponding Theorem 4.16 in \cite{FHS1} one realizes that the only thing one has to prove is that each point of $\partial_i^+ M $ is deformably timelikely accessible. But that is an easy exercise using causal convexity of $\overline{i(M)}$. \hfill \qed

\bigskip
\bigskip

\v{\bf Acknowledgements.} It is my pleasure to thank for many helpful comments about the subjects of this article by Christian B\"ar, Claudio Dappiaggi, Jos\'e Luis Flores, Steven Harris and Miguel S\'anchez.


\begin{thebibliography}{99}

\bibitem{AC}
Michael T. Anderson, Piotr T. Chrusciel: {\em Asymptotically simple solutions of the vacuum Einstein equations in even dimensions}. Commun.Math.Phys.260:557-577 (2005). arXiv: gr-qc/0412020

\bibitem{BEE}
John K. Beem, Paul Ehrlich, Kevin Easley: {\em Global Lorentzian geometry}, 2nd edition. CRC Press (1996)

\bibitem{BS3}
Antonio N. Bernal, Miguel S\'anchez: {\em Further results on the smoothability of Cauchy hypersurfaces and Cauchy time functions}, Letters in Mathematical Physics 77 (2), 183-197 (2006). arXiv: /gr-qc/0512095

\bibitem{CC}
Y.~Choquet-Bruhat, D.~Christodoulou, \emph{Existence of global solutions of the Yang-Mills, Higgs and spinor field equations in 3+1 dimensions}, Ann. Sci. \'Ecole Norm. Sup. (4) \textbf{14}, no. 4, 481--506 (1981)

\bibitem{pC}
Piotr T. Chrusciel: {\em Elements of Causality theory}, Preprint University of Viena (UWThPh-2011-32), arXiv: 1110.6706 (2011)

\bibitem{CGM}
Piotr T. Chrusciel, James D. E. Grant, Ettore Minguzzi: {\em On differentiability of volume time functions}, preprint, arxiv: 1301.2909

\bibitem{D}
Sergio Dain: {\em Initial data for stationary space-times near space-like infinity}, Class.
Quantum Grav. 18 (2001), 4329–4338, arXiv: gr-qc/0107018.

\bibitem{DHP}
Claudio Dappiaggi, Thomas Hack, Nicola Pinamonti: {\em Approximate KMS states for scalar and spinor fields in Friedmann–Robertson–Walker spacetimes}, Annales Henri Poincare 12 (8), 1449-1489

\bibitem{DS}
Thibault Damour, Berndt Schmidt: {\em Reliability of perturbation theory in general relativity},
Jour. Math. Phys. 31 (1990), 2441–2453.

\bibitem{FH}
Jos\'e Luis Flores, Steven G. Harris: {\em Topology of the causal boundary for standard static spacetimes}, Class.Quant.Grav.24:1211-1260 (2007). arXiv: gr-qc/0607049

\bibitem{FHS1}
Jos\'e Luis Flores, Jonatan Herrera, Miguel S\'anchez: {\em On the final definition of the causal boundary and its relation with the conformal boundary}, Adv. Theor. Math. Phys. Volume 15, Issue 4, 991-1058 (2011). arXiv:1001.3270

\bibitem{FHS2}
Jos\'e Luis Flores, Jonatan Herrera, Miguel S\'anchez: {\em Gromov, Cauchy and causal boundaries for Riemannian, Finslerian and Lorentzian manifolds}. Memoirs Amer. Mat. Soc. 226, No. 1064 (2013). arXiv: 1011.1154


\bibitem{cF}
Charles Frances: {\em About geometrically maximal manifolds}, Journal of Topology. Vol 5. no 2., 293-322 (2012)

\bibitem{hF}
Helmut Friedrich: {\em On the existence of $n$-geodesically complete or future complete solutions of Einstein’s field equations with smooth asymptotic structure}, Commun. Math. Phys. 107, 587–609 (1986)


\bibitem{G}
Robert Geroch: {\em Domain of dependence}. J. Math. Phys. 11, 437 (1970)

\bibitem{GKP}
Robert Geroch, E.H. Kronheimer, Roger Penrose: {\em Ideal points in space-time}. Proc. R. Soc. Lond. A vol. 327 no. 157, pp. 545-567 (1972)

\bibitem{GM}
Nicolas Ginoux, Olaf M\"uller: {\em Massless Dirac-Maxwell systems in asymptotically flat spacetimes}, preprint, arXiv: 1407.1177 (2014)

\bibitem{H}
Steven G. Harris: {\em Causal Boundary for Standard Static Spacetimes}, Nonlinear Anal., 47, 2971 - 2981 (2001)

\bibitem{mH}
Marc Herzlich: {\em Compactification conforme des variet\'es asymptotiquement plats}, Bull. Soc. Math. France 125, pp. 55-92 (1997)

\bibitem{JS}
Miguel Angel Javaloyes, Miguel S\'anchez: {\em Existence of standard splitting for conformally stationary spacetimes}, Class. Quant. Grav. 25, 168001 (2008). arXiv: 0806.0812

\bibitem{LR}
Hans Lindblad, Igor Rodnianski: {\em The global stability of the Minkowski space-time in harmonic gauge}. Annals of Math 171, no. 3, 1401-1477 (2010). arXiv: math/0411109

\bibitem{wML}
William McLean: {\em Strongly Elliptic Systems and Boundary Integral Equations}, Cambridge University Press (2000)

\bibitem{vM}
Valter Moretti: {\em Uniqueness theorem for BMS-invariant states of scalar
QFT on the null boundary of asymptotically flat space-times and bulk-boundary observable algebra correspondence}, Commun.Math.Phys. 268, 727-756 (2006). arXiv.org: gr-qc/0512049

\bibitem{oM3}
Olaf M\"uller: {\em Horizons}, preprint, arxiv: 1111.4571 (2012)

\bibitem{PR}
Roger Penrose, Wolfgang Rindler: {\em Spinors and Spacetime vol.2: Spinor and Twistor Methods in Space-Time Geometry}. Cambridge University Press (1986)

\bibitem{SB}
Walter Simon, Robert Beig: {\em The multipole structure of stationary space–times}, Jour. Math. Phys. 24, 1163–1171 (1983)



\end{thebibliography}
\end{document}